\documentclass[11pt]{article}%
\usepackage{amsmath}
\usepackage{amsfonts}
\usepackage{amssymb}
\usepackage{graphicx}
\usepackage{a4wide}%
\providecommand{\U}[1]{\protect\rule{.1in}{.1in}}
\makeatletter
\renewcommand\@biblabel[1]{#1.}
\makeatother

\begin{document}

\title{A critical view on invexity}
\author{Constantin Z\u{a}linescu\thanks{University ``Al.I.Cuza'' Ia\c{s}i, Faculty of
Mathematics, and Institute of Mathematics Octav Mayer, Ia\c{s}i, Romania;
email: \texttt{zalinesc@uaic.ro} }}
\date{}
\maketitle

\noindent\textbf{Abstract:} {The aim of this note is to emphasize the fact
that in many papers on invexity published in prestigious journals there are
not clear definitions, trivial or not clear statements and wrong proofs. We
also point out the unprofessional way of answering readers' questions by some
authors. We think that this is caused mainly by the lack of criticism of the
invexity community
}


\section{Introduction}

There are a lot of articles written on invexity. A search in MathScieNet gives
350 articles having in their titles the words ``invex, invexity, preinvex"
(310 in Zentralblatt f\"{u}r Mathematik).

In my opinion, the number of articles devoted to invexity (invex and
generalized invex functions) is too big with respect to its importance. In
fact, when I read the first time about (classic!?) invex functions, in the
differentiable case, I realized that this is another way of saying that any
critical (or stationary) point of the function, that is, a point at which the
differential is $0$, is necessarily a global minimum. Then the notion was
generalized to non differentiable functions, but saying the same thing: every
critical point (in the sense that $0$ is in a certain type of subdifferential
at that point) is a global minimum. If one looks to the applications of the
results in Ref.~1
we see that they state something like: every local solution is a global one.

Rapidly one had generalizations of the notion: quasiinvex functions, preinvex
functions, and so on. The common feature for many papers on invexity and its
generalizations is the lack of clarity of the notions and statements of the
results and doubtful proofs; when the proofs are correct many of them are
trivial. After the first draft of this note was written (and after the change
of messages with some authors of the cited papers) I had a closer look to
several reviews in Mathematical Reviews and Zentralblatt f\"{u}r Mathematik
(referred in Section 4), some of them related to the quoted articles; it seems
that the opinions of the reviewers didn't influence the authors of papers on invexity.

\section{About statements and proofs}

Let us consider the following text quoted from Ref.~1
(Ref.~1
appears as being cited 37 times in Google Scholar and 4 times in MathSciNet at
the time of writing this note):

\medskip

\textquotedblleft\textbf{Definition 1.1}. See Refs. 1--2. A set $K\subseteq
\mathbb{R}^{n}$ is said to be invex if there exists a vector function
$\eta:\mathbb{R}^{n}\times\mathbb{R}^{n}\rightarrow\mathbb{R}^{n}$ such that%
\[
x,y\in K,\ \lambda\in\lbrack0,1]\Rightarrow y+\lambda\eta(x,y)\in K.
\]

\textbf{Remark 1.1.} A convex set is an invex set; i.e., take $\eta(x,y)=x-y$.
But the converse does not hold.\textquotedblright\footnote{To see the precise
coordinates of the references referred in the quoted texts one might consult
the reviews on MathSciNet of the corresponding articles listed at the end of
this note.}

\medskip

Note on this remark: Of course, the converse does not hold because by
Definition 1.1 in Ref.~1
quoted above, any set is invex: just take $\eta(x,y):=0$. The honest
definition is: let $\eta:\mathbb{R}^{n}\times\mathbb{R}^{n}\rightarrow
\mathbb{R}^{n}$ be a function. The set $K\subset\mathbb{R}^{n}$ is said to be
$\eta$-invex if ...

\medskip

Or there is another formulation: One says that $f$ is $(\theta,\alpha)$-$d$
invex if there exist $\theta$ etc such that $f$ satisfies a certain condition
involving $\theta$ etc. Of course, correctly is to introduce first
$\theta,\alpha,d$ and after that to say that $f$ is $(\theta,\alpha)$-$d$
invex if ...

Let us quote from Ref.~2
:

\medskip

\textquotedblleft\textbf{Theorem 8.} A function $f:R^{n}\rightarrow R$ is
$B$-$(0,r)$-invex ($B$-$(0,r)$-incave) with respect to $\eta$ and $b$ on
$R^{n}$ if and only if its every stationary point is a global minimum
(maximum) in $R^{n}$.\textquotedblright

\medskip

At least two remarks are in order with respect to this statement. The first
one: If the statement is true, $f$ is $B$-$(0,r)$-invex with respect to $\eta$
and $b$ (in the sense of Definition 1 in Ref.~2
if and only if $f$ is invex (because, as seen above, the invexity of the
differentiable function $f:\mathbb{R}^{n}\rightarrow\mathbb{R}$ is equivalent
to the fact that every stationary point is a global minimum); so which is the
need to introduce $B$-$(0,r)$-invexity? The second one: The statement gives
the impression that the functions $\eta$ and $b$, as well as $r\in\mathbb{R}$,
are given. Consider $n=1$, $r=0$, $b(x,u):=1$ and $\eta(x,u):=u^{-1}\left(
x^{2}-u^{2}\right)  +\operatorname*{sgn}u$ for $u\neq0$, $\eta(x,u):=0$ for
$u=0$. Taking $f(x)=\tfrac{1}{2}x^{2}$ for $x\in\mathbb{R}$ we see that every
stationary point of $f$ (that is, $u=0$) is a global minimum, but $f$ is not
$B$-$(0,r)$-invex with respect to $\eta$ and $b$ on $\mathbb{R}$. Let us quote
also from Definition 8 in Ref.~3
:

\textquotedblleft Let $S\subset R^{n}$ be a nonempty invex set with respect to
$\eta$. A function $f:S\rightarrow R$ is said to be pre-invex with respect to
$\eta$ if, there exists a vector-valued function $\eta:S\times S\rightarrow
R^{n}$ such that the relation ...\textquotedblright;

Definition 9 in Ref.~3
is obtained by changing pre-invex by invex. So, first $\eta$ is given, and one
line below one asks the existence of an $\eta$.

I am not used with this kind of text in mathematics. If invexity is not a
domain of mathematics it is advisable to say it explicitly. Why is this
important? Because we are judged in comparison with other mathematicians for
getting jobs, for promotions, for getting grants. Even if one declares that
invexity is not a domain of mathematics, this does not change a lot the
situation because now one asks for interdisciplinarity.

Why did I ask if invexity is a topic in mathematics? Because I had the
impression that the next text quoted from Ref.~4
is not a mathematical text (Ref.~4
appears as being cited 21 times in Google Scholar and 3 times in MathSciNet):

\medskip

\textquotedblleft\textbf{Remark 2.3.} We will show that Assumption C holds if

$\eta(x,y)=x-y+o(x-y)$. In fact, the following two equalities hold:

(i) $\eta(y,y+\lambda\eta(x,y))$

$=\eta(y,y+\lambda(x-y+o(x-y))$

$=-\lambda(x-y+o(x-y))+o(\lambda(x-y+o(x-y)))$

$=-\lambda\lbrack x-y+o(x-y)+o(x-y+o(x-y))]$

$=-\lambda\lbrack x-y+o(x-y)]$

$=-\lambda\eta(x,y);$" \footnote{Related to this piece of non-mathematics let
us quote from what the authors say in the first footnote on the first page of
Ref.~4
: \textquotedblleft The authors are thankful to ..., and \textbf{three
anonymous referees} for their many valuable comments on an early version of
this paper. The authors are also grateful to Professor B.D. Craven for some
discussion on this paper\textquotedblright. (C.Z.' emphasis.)}

\medskip

\noindent(I didn't quote the second equality which is of the same type).

Seeing this I sent messages to the authors of Ref.~4 asking:

\medskip

\textquotedblleft What do you mean by $\eta(x,y)=x-y+o\left(  \left\Vert
x-y\right\Vert \right)  $ in Remark 2.3?\textquotedblright,

\medskip

\noindent but the answers were unsatisfactory. I don't cite the
answers here being contained in private correspondence.

\medskip

Which is this Condition C (or Assumption C)?

I quote again from page 610 of Ref.~4
(see also page 116 of Ref.~5):

\medskip

\textquotedblleft\textbf{Assumption C.} See Ref.6. Let $\eta:X\times
X\rightarrow\mathbb{R}^{n}$. Then, for any $x,y\in\mathbb{R}^{n}$

and for any $\lambda\in\lbrack0,1]$,

$\eta(y,y+\lambda\eta(x,y))=-\lambda\eta(x,y)$, $\eta(x,y+\lambda
\eta(x,y))=(1-\lambda)\eta(x,y).$\textquotedblright

\medskip

Somewhere it is written that $X\subset\mathbb{R}^{n}$; probably for the
authors it is not very important to speak about $\eta(x,y)$ when $x$ or $y$ is
not in $X$. Let us consider $X=\mathbb{R}^{n}$. (Also note that in Definition
2.4 of Ref.~4
$\eta:X\times X\rightarrow\mathbb{R}$, that is, $\eta$ takes its values in
$\mathbb{R}$ instead of $\mathbb{R}^{n}.$)

To see that the condition $\eta(x,y)=x-y+o(x-y)$ does not imply
Condition C let us consider
$\eta:\mathbb{R}\times\mathbb{R}\rightarrow\mathbb{R}$ be defined by
$\eta(x,y)=x-y+(x-y)^{2}$. Of course $\eta$ satisfies the hypothesis
of the statement in Remark 2.3 of Ref.~4 (above). For $\lambda=1$,
the first relation of Condition C is equivalent to each of the next
ones: $-\eta(x,y)+\left(  \eta(x,y)\right)  ^{2}=-\eta(x,y)$,
$\left( \eta(x,y)\right)  ^{2}=0$, $\eta(x,y)=0$, $(x-y)(1+x-y)=0$,
$x-y\in\{0,-1\}$. So, taking $x=0$, $y=2$ and $\lambda=1$ one sees
that $\eta$ does not verify Condition C.

In Ref.~5
one finds:

\medskip

\textquotedblleft\textbf{Example 2.1}. Let

$\eta(x,y)=\left\{
\begin{array}
[c]{ccc}%
x-y & \text{if} & x\geq0,\ y\geq0;\\
x-y & \text{if} & x\leq0,\ y\leq0;\\
-2-y & \text{if} & x>0,\ y\leq0;\\
2-y & \text{if} & x\leq0,\ y>0.
\end{array}
\right.  $

Then, it is easy to verify that $\eta$ satisfies Condition
C.\textquotedblright

\medskip

A similar example is the following quoted from Ref.~4
(which is very close to that quoted above from Ref.~5, as well as to
Example 1.1 in Ref.~1):

\medskip

\textquotedblleft\textbf{Example 2.2}. Let

$f(x)=-\left\vert x\right\vert ,\ \forall x\in K=[-2,2]$, and let

$\eta(x,y)=\left\{
\begin{array}
[c]{ccc}%
x-y & \text{if} & x\geq0,\ y\geq0;\\
x-y & \text{if} & x<0,\ y<0;\\
-2-y & \text{if} & x>0,\ y\leq0;\\
2-y & \text{if} & x\leq0,\ y>0.
\end{array}
\right.  $

Then, it is easy to verify that $f$ is invex with respect to $\eta$ on $K$ and
that $f$ and $\eta$ satisfy Assumptions A and C. However, $f$ is not
convex.\textquotedblright

\medskip

The authors seem to not realize that $\eta$ defined in these two
examples is not a function because $\eta(2,0)$ gives $2$ using the
first expression and $-2$ using the third expression. A possible
modification for $\eta$ defined in Example 2.1 of Ref.~5 could be:
\[
\eta(x,y)=\left\{
\begin{array}
[c]{ccc}%
x-y & \text{if} & xy\geq0,\\
2-y & \text{if} & xy<0.
\end{array}
\right.
\]

Take $x,y\in\mathbb{R}$ with $x>0>y$ and $\lambda\in\lbrack0,1]$ and let us
look to the second relation in Condition C. We have that $\eta(x,y)=2-y$ and
$y^{\prime}:=y+\lambda\eta(x,y)=y+\lambda(2-y)$. Assuming that $y^{\prime}%
\geq0$, then $\eta(x,y^{\prime})=x-y^{\prime}=x-2+(1-\lambda)\eta(x,y)$.
Hence, in such a situation ($y<0$ and $y+\lambda(2-y)\geq0$) one has
$\eta(x,y+\lambda\eta(x,y))=(1-\lambda)\eta(x,y)$ if and only if $x=2$. Is it
possible to have $y<0$, $y+\lambda(2-y)\geq0$ and $\lambda\in\lbrack0,1]?$ The
answer is YES! Just take $\lambda=1$. Hence, for $x=1$, $y=-1$ and $\lambda=1$
the second relation in Condition C is not verified because $\eta(1,-1)=3$,
$\eta(1,-1+3)=\eta(1,2)=-1\neq0$.

In fact an adequate modification of the function $\eta$ in Example 2.1 from
Ref.~5
(or Example 2.2 in Ref.~4) is
\[
\eta(x,y)=\left\{
\begin{array}
[c]{ccl}%
x-y & \text{if} & xy\geq0,\\
2-y & \text{if} & x<0,\ y>0,\\
-2-y & \text{if} & x>0,\ y<0.
\end{array}
\right.
\]
The function $\eta$ defined in this way satisfies indeed Condition C.

Somewhere (say [S]) it was said that
\begin{equation}
\eta(y+\lambda_{2}\eta(x,y),y+\lambda_{1}\eta(x,y))=(\lambda_{2}-\lambda
_{1})\eta(x,y)\quad\forall x,y\in\mathbb{R}^{n},\ \forall\lambda_{1}%
,\lambda_{2}\in\lbrack0,1], \label{r1}%
\end{equation}
whenever $\eta$ verifies Condition C; and for this the proof of Theorem 3.1 in
Ref.~4
was cited. In fact I was determined by [S] to look at Ref.~4
and Ref.~5.
Of course, relation (\ref{r1}) is nice and good to have; moreover,
for $\lambda_{2}=0$ one recovers the first relation in Condition C.

Looking at the proof of Theorem 3.1 in Ref.~4
(but one can look also at the proof of Theorem 2.1 in Ref.~4
for the same text), one observes that one takes
$0<\lambda_{2}<\lambda_{1}<1$ and one obtains relation (14) of Ref.\
4 I am quoting below:

\medskip

``$\eta(y+\lambda_{1}\eta(x,y),y+\lambda_{2}\eta(x,y))$

$=\eta(y+\lambda_{1}\eta(x,y),y+\lambda_{1}\eta(x,y)-(\lambda_{1}-\lambda
_{2})\eta(x,y))$

$=\eta(y+\lambda_{1}\eta(x,y),y+\lambda_{1}\eta(x,y)+\eta(y,y+(\lambda
_{1}-\lambda_{2})\eta(x,y)))$

$=-\eta(y,y+(\lambda_{1}-\lambda_{2}\eta(x,y))$

$=(\lambda_{1}-\lambda_{2})\eta(x,y).$''

\medskip

The first equality is obvious, the second as well as the fourth follow from
the first relation in Condition C [however it is $=-\eta(y,y+(\lambda
_{1}-\lambda_{2})\eta(x,y))$ instead of $=-\eta(y,y+(\lambda_{1}-\lambda
_{2}\eta(x,y))$]. What is used to obtain the third equality? Setting
$y^{\prime}:=y+\lambda_{1}\eta(x,y)$, the expression on the third line becomes
$\eta(y^{\prime},y^{\prime}+\eta(y,y+(\lambda_{1}-\lambda_{2})\eta(x,y)))$. In
order to get the expression on the fourth line (using directly Condition C) we
should have $\eta(y^{\prime\prime},y^{\prime\prime}+\eta(y,y^{\prime\prime}))$
with $y^{\prime\prime}:=y+(\lambda_{1}-\lambda_{2})\eta(x,y)$. Is $y^{\prime
}=y^{\prime\prime}?$ In fact $y^{\prime}=y^{\prime\prime}$ if and only if
$\eta(x,y)=0$ or $\lambda_{2}=0$. \medskip

Maybe (\ref{r1}) is true whenever Condition C holds, but some additional
arguments must be provided.

I do not propose myself to mention all doubtful sentences or statements in
articles about invexity, but the majority I had occasion to browse are like that.

\section{About the triviality of results and generalizations}

Another problem with invexity is given by the triviality of some results or
generalizations. Let us mention some of them found in recent articles
published in prestigious journals.

It is well known that for a G\^{a}teaux differentiable function
$f:D\rightarrow\mathbb{R}$ with $D$ an open subset of a normed
vector space (but $X$ could be a topological vector space), for any
(distinct) points $a,b\in D$ with $[a,b]:=\{\lambda
a+(1-\lambda)b\mid\lambda\in\lbrack 0,1]\}\subset D$, there exists
$c\in{}]a,b[{}:=\{\lambda a+(1-\lambda )b\mid\lambda\in(0,1)\}$ such
that $f(b)-f(a)=\nabla f(c)(b-a)$ (the proof being immediate using
the real-valued function $\varphi$ defined by $\varphi(t):=f\left(
ta+(1-t)b\right)  $ for those $t\in\mathbb{R}$ with $ta+(1-t)b\in
D$). Which are the main results obtained in Ref.~3?
I don't speak about Theorems 11 and 12 which are just rewriting of
the definitions of convexity and pre-invexity (in Theorem 11 of
Ref.~3
no need of differentiability of $f$). Let us quote Theorem 14 in Ref.~3:

\medskip

``\textbf{Theorem 14.} Let $S\subset R^{n}$ be a nonempty invex set with
respect to $\eta:S\times S\rightarrow R^{n}$, and $P_{ab}$ be an arbitrary
$\eta$-path contained in $\operatorname*{int}S$. Moreover, we assume that
$f:S\rightarrow R$ is defined on $S$ and differentiable on
$\operatorname*{int}S$. Then, for any $a,b\in S$, there exists $c\in
P_{ab}^{0}$ such that the following relation

$f(a+\eta(b,a))-f(a)=[\eta(b,a)]^{T}\nabla f(c)\qquad(8)$

\noindent holds.\textquotedblright\

\medskip

Because (c.f.\ Definition 5 in Ref.~3)
$P_{ab}:=[a,a+\eta(a,b)]$ and $P_{ab}^{0}:={}]a,a+\eta(a,b)[$, we
see that
Theorem 14 in Ref.~3
 is an immediate consequence of the usual mean-value theorem mentioned above.
(Note that it is not said what kind of differentiability is asked for $f$ ---
G\^{a}teaux or Fr\'{e}chet.) Probably the next paper will deal with such a
result in infinite dimensional spaces, then with $\alpha$-$\eta$ invex
functions (mentioned below). Theorem 17 in Ref.~3
deals with a Taylor's expansion (of order 2) for $f$. Other ``{important}''
results (Theorems 21, 22) are immediate consequences of known results for
derivable functions of one real variable. They could constitute easy exercises
for students following a first course in analysis.

Another example in this sense is provided by Ref.~6.
As seen in the title of Ref.~6,
there is some $G$ there. What is it? It is a function defined on a
certain set $A\subset\mathbb{R}$ with values in $\mathbb{R}$ which
is increasing ($s,t\in A$, $s<t$ $\Rightarrow$ $G(s)<G(t)$), and
moreover $G$ is differentiable. In fact $G$ is defined on the image
of a real-valued function $f$ defined in its turn on an $\eta$-invex
set $X\subset\mathbb{R}^{n}$. (By the way, if $A=I_{f}(X)$ is
$\{0,1\}$, what does differentiability of $G$ mean?) One defines
$G$-invex and $G$-pre-invex functions. Let us quote Definition 3 in
Ref.~6:

\medskip

\textquotedblleft\textbf{Definition 3.} Let $X$ be a nonempty invex (with
respect to $\eta$) subset of $R^{n}$ and $f:X\rightarrow R$ be a
differentiable function defined on $X$. Further, we assume that there exists a
differentiable real-valued increasing function $G:I_{f}(X)\rightarrow R$. Then
$f$ is said to be (strictly) $G$-invex at $u\in X$ on $X$ with respect to
$\eta$ if there exists a vector-valued function $\eta:X\times X\rightarrow
R^{n}$ such that, for all $x\in X$ $(x\neq u)$, $G(f(x))-G(f(u))\geq
G^{\prime}(f(u))\nabla f(u)\eta(x,u)\qquad(>).\qquad(2)$

If (2) is satisfied for any $u\in X$ then $f$ is $G$-invex on $X$ with respect
to $\eta$.\textquotedblright

\medskip

Taking into account that for $f$ (Fr\'{e}chet) differentiable one has
$\nabla h(u)=G^{\prime}(f(u))\nabla f(u)$, where $h:=G\circ f$, the
inequality above says that $h(x)-h(u)\geq\nabla h(u)\eta(x,u)$.
Having this inequality for all $x,u\in X$ this means that $h$ is
invex. So, one can say simply that $f$ is $G$-invex (at $u$) if
$G\circ f$ is invex (at $u$). This simple remark is not made in
Ref.~6,
 but one has (quoted from Ref.~6
):

\medskip

\textquotedblleft We remark that the $G$-invexity assumption
generalizes a hypothesis of Avriel et al.\ [6], Avriel [7], Hanson
[11] and Antczak [3] for differentiable functions. Thus, the
following remarks are true:

\textbf{Remark 5.} In the case when $\eta(x,u)=x-u$, we obtain a definition of
a differentiable $G$-convex function introduced Avriel et al. [6].

\textbf{Remark 6.} Every invex function with respect to $\eta$ introduced by
Hanson [11] is $G$-invex with respect to the same function $\eta$, where
$G:I_{f}(X)\rightarrow R$ is defined by $G(a)\equiv a$. The converse result
is, in general, not true (see also Remark 13 and Example 14).

\textbf{Remark 7.} Every $r$-invex function with respect to $\eta$ introduced
by Antczak [1,3] is $G$-invex with respect to the same function $\eta$, where
$G:I_{f}(X)\rightarrow R$ is defined by $G(a)=e^{ra}$, where $r$ is any finite
real number.\textquotedblright\

\medskip

(However, note that for $r\leq0$ the function $G$ defined by $G(a)=e^{ra}$ is
not increasing.)

It is suggestive to quote also the definition a $G$-pre-invex function (but
probably the reader already guesses it):

\medskip

\textquotedblleft\textbf{Definition 9.} Let $X$ be a nonempty invex (with
respect to $\eta$) subset of $R^{n}$. A function $f:X\rightarrow R$ is said to
be (strictly) $G$-pre-invex at $u$ on $X$ with respect to $\eta$ if there
exist a continuous real-valued increasing function $G:I_{f}(X)\rightarrow R$
and a vector-valued function $\eta:X\times X\rightarrow R^{n}$ such that for
all $x\in X$ $(x\neq u)$,

$f(u+\lambda\eta(x,u))\leq G^{-1}(\lambda G(f(x))+(1-\lambda)G(f(u)))\qquad
(<).\qquad(3)$

If (2) is satisfied for any $u\in X$ then $f$ is $G$-pre-invex on $X$ with
respect to $\eta$.\textquotedblright

\medskip

Of course, the author does not (want to) observe that this means that
$h:=G\circ f$ is pre-invex at $u$. What does he obtain in Theorem 10 of Ref.~6
? He obtains that $f$ is $G$-invex provided $f$ and $G$ are differentiable and
$f$ is $G$-pre-invex. I quote from page 646 in Ref.~1
:

\medskip\textquotedblleft Recently, Pini (Ref.~6) showed that, if $f$ is
defined on an invex set $K\subseteq\mathbb{R}^{n}$ and if it is preinvex and
differentiable, then $f$ is also invex with respect to $\eta$; i.e.,
$f(y)-f(x)\geq\eta(y,x)^{T}\nabla f(x)$\textquotedblright.

\medskip Of course, in Ref.~6
one gives a detailed proof. At page 646 of Ref.~1
one continues with:

\medskip\textquotedblleft But the converse is not true in general. A
counterexample was given in Ref.~6. However, Mohan and Neogy (Ref.\
9) proved that a differentiable invex function is also preinvex
under the following condition. \textbf{Condition C.}
...\textquotedblright\

\medskip At page 1620 of Ref.~6
one says:

\medskip\textquotedblleft The converse result is not true in general, that is,
there exist $G$-invex functions with respect to $\eta$ which are not
$G$-pre-invex with respect to the same function $\eta$. To prove the converse
theorem the function $\eta$ should satisfy the following condition C (see
[16]). \textbf{Condition C}. ...\textquotedblright\

\medskip Of course one states Theorem 11 and one gives a detailed proof. As a
conclusion for paper Ref.~6
: The function $f$ is $G$-``word" if $G\circ f$ is ``word". If an existing
result holds for ``word" then in Ref.~6
one has a result for $G$-``word" with detailed proof. And this is published in
a prestigious journal.

\medskip

The case of Refs.~3,~6
is not singular. Let us have a look to Ref.~7
and its follower Ref.~8
. Let us quote first from Ref.~7
two interesting phrases:

\medskip\textquotedblleft In recent years, the concept of convexity has been
generalized and extended in several directions using \textbf{novel and
innovative} techniques\textquotedblright\

\medskip\noindent and

\medskip\textquotedblleft Motivated and inspired by the research going on in
this \textbf{fascinating field}, we introduce a new class of generalized
functions\textquotedblright.

\medskip Let us quote again from page 698 of Ref.~7
:

\medskip

\textquotedblleft Let $K$ be a nonempty closed set in a real Hilbert
space $H$. We denote by $\left\langle .,.\right\rangle $ and
$\left\Vert .\right\Vert $ the inner product and norm respectively.
Let $F:K\rightarrow H$ and $\eta(.,.):K\times K\rightarrow R$ be
continuous functions. Let $\alpha:K\times K\rightarrow
R\setminus\{0\}$ be a bifunction. First of all, we recall the
following well-known results and concepts.

\textbf{Definition 2.1.} Let $u\in K$. Then the set $K$ is said to be $\alpha
$-invex at $u$ with respect to $\eta(.,.)$ and $\alpha(.,.)$, if, for all
$u,v\in K,t\in\lbrack0,1]$, $u+t\alpha(v,u)\eta(v,u)\in K$. $K$ is said to be
an $\alpha$-invex set with respect to $\eta$ and $\alpha$, if $K$ is $\alpha
$-invex at each $u\in K$. The $\alpha$-invex set $K$ is also called
$\alpha\eta$-connected set. Note that the convex set with $\alpha(v,u)=1$ and
$\eta(v,u)=v-u$ is an invex set, but the converse is not
true.\textquotedblright

\medskip

First note that $u+t\alpha(v,u)\eta(v,u)$ above does not make sense if
$H\neq\mathbb{R}$ because $u\in H$ and $t\alpha(v,u)\eta(v,u)\in\mathbb{R}$;
next, if $\eta(.,.):K\times K\rightarrow H$ (as in Ref.~8
, then $K$ is an $\alpha$-invex set with respect to $\eta$ iff $K$ is
$\eta^{\prime}$-invex, where $\eta^{\prime}:=\alpha\eta$ (apparently not
observed in Refs.~7, 8). Of course, in Definition 2.2 of Ref.~7
one says:

\medskip\textquotedblleft The function $F$ on the $\alpha$-invex set $K$ is
said to be $\alpha$-preinvex with respect to $\alpha$ and $\eta$, if
$F(u+t\alpha(v,u)\eta(v,u))\leq(1-t)F(u)+tF(v),\ \forall u,v\in K,\ t\in
\lbrack0,1]$\textquotedblright,

\medskip\noindent that is (I say), $F$ is $\eta^{\prime}$-preinvex (however,
one must take $F:K\rightarrow R$ as in Ref.~8
instead of $F:K\rightarrow H)$. In a similar way one obtains the corresponding
definitions for \textquotedblleft$\alpha$-invex\textquotedblright\ replaced by
\textquotedblleft quasi $\alpha$-preinvex\textquotedblright\ (see Definition
2.3 in Ref.~7
, \textquotedblleft logarithmic $\alpha$-preinvex\textquotedblright\ (see
Definition 2.4 in Ref.~7
), \textquotedblleft pseudo $\alpha$-preinvex\textquotedblright\ (see
Definition 2.5 in Ref.~7
) from the definitions without $\alpha$. (Note the interesting inequality
$\max\{F(u),F(v)\}<\max\{F(u),F(v)\}$ from the displayed relation after
Definition 2.4 in Ref.~7
.) Maybe the next one is an exception:

\medskip

\textquotedblleft\textbf{Definition 2.6.} The differentiable function $F$ on
$K$ is said to be an $\alpha$-invex function with respect to $\alpha$ and
$\eta$, if

$F(v)-F(u)\geq\left\langle \alpha(v,u)F^{\prime}(u),\eta(v,u)\right\rangle
,\qquad\forall u,v\in K$, \noindent where $F^{\prime}(u)$ is the differential
of $F$ at $u\in K$. \textbf{The concepts of the }$\alpha$\textbf{-invex and
}$\alpha$\textbf{-preinvex functions have played very important role in the
development of convex programming; see [6,7]}. Note that for $\alpha(v,u)=1$,
Definition 2.6 is mainly due to Hanson [1]\textquotedblright.

\medskip

Unfortunately not, even in this case, $F$ is $\alpha$-invex with respect to
$\alpha$ and $\eta$ iff $F$ is $\eta^{\prime}$-invex. What is new and
surprising for me is the emphasized text above.

Similar remarks are valid for the notions of \textquotedblleft$\alpha\eta
$-monotone\textquotedblright, \textquotedblleft strictly $\alpha\eta
$-monotone\textquotedblright, \textquotedblleft$\alpha\eta$%
-pseudomonotone\textquotedblright, \textquotedblleft quasi $\alpha\eta
$-monotone\textquotedblright, \textquotedblleft strictly $\alpha\eta
$-pseudomonotone\textquotedblright\ referred to an operator $T:K\rightarrow H$
(defined in Definition 2.7 in Ref.~7
).

However, there are some notions which do not correspond to those for
$\eta^{\prime}:=\alpha\eta$. These are those containing the word
\textquotedblleft strongly\textquotedblright\ in their definition:
\textquotedblleft strongly $\alpha\eta$-monotone\textquotedblright\ and
\textquotedblleft strongly $\alpha\eta$-pseudomonotone\textquotedblright%
\ operators (see Definition 2.7 in Ref.~7
) as well as \textquotedblleft strongly $\alpha$-preinvex\textquotedblright%
\ (see Definition 2.8 in Ref.~7
), \textquotedblleft strongly $\alpha$-invex\textquotedblright\ (see
Definition 2.9 in Ref.~7
), \textquotedblleft strongly pseudo $\alpha\eta$-invex\textquotedblright%
\ (see Definition 2.10 in Ref.~7
) and \textquotedblleft strongly quasi $\alpha$-invex\textquotedblright\ (see
Definition 2.11 in Ref.~7
) functions. The results which refer to these notions are Theorems 3.1--3.5 in
Ref.~7
. I do not propose myself to verify the correctness of these results
(however see Example 6.1 in Ref.~8
, but some of them probably are not true having in mind that Theorems 6.1 and
6.4 in Ref.~8 
give alternative formulations for the sufficiency parts of Theorems 3.2 and
3.5 in Ref.~7
, respectively. What I want to point out are the following facts:

1) If $\eta(u,u)=0$ for some $u\in K$ then there do not exist pseudo $\alpha
$-preinvex, strictly $\alpha$-invex and strictly pseudo $\alpha$-invex
functions with respect to $\alpha$ and $\eta$, as well as strictly $\alpha
\eta$-monotone and strictly pseudo $\alpha\eta$-monotone operators. Note that
if $\alpha$ and $\eta$ satisfy Condition C at page 702 of Ref.~7
 or condition (ii) in Theorem 6.1 of Ref.~8
 then $\eta(u,u)=0$ for every $u\in K$.

2) In some proofs of the statements in Refs.~7 and 8 
one uses the relation $g(1)-g(0)=\int_{0}^{1}g^{\prime}(t)dt$, where $g$ is a
real-valued derivable function on a subset of $\mathbb{R}$ containing $[0,1]$.
In fact $g(t):=F(u+t\alpha(v,u)\eta(v,u))$ for $t\in\lbrack0,1]$, where $F$ is
differentiable. It is a well-known fact that the formula $g(1)-g(0)=\int
_{0}^{1}g^{\prime}(t)dt$ might not be true if $g^{\prime}$ is not Riemann
integrable on $[0,1]$. As an example take $g(t):=t^{2}\sin t^{-2}$ for $t\in
{}]0,1]$, $g(0):=0$.

3) In Theorems 6.1--6.4 of Ref.~8
one uses the condition \textquotedblleft$\alpha(u,u+t\alpha(v,u)\eta
(v,u))=t\alpha(v,u),\ \forall u,v\in K,\ t\in\lbrack0
,1]$\textquotedblright.
Taking $t=0$, this implies that $\alpha(u,u)=0$ for every $u\in K$,
contradicting the assumption made before Definition 1.1 in Ref.~8
that $\alpha$ takes its values in $R\setminus\{0\}$. This shows that the
domain of applicability of Theorems 6.1--6.4 in Ref.~8
is the empty set.

\section{Conclusions}

In this note we pointed out that several papers published in
prestigious journals contain important drawbacks in the formulation
of the notions and in the statements of the results, as well as very
serious mistakes in the proofs. Also, there are many trivial
generalizations of notions and results. In this sense it is useful
to mention that there are several reviews in Mathematical Reviews
and Zentralblatt f\"{u}r Mathematik which are concordant with our
opinions; let us cite the reviews \emph{MR1989930 (2004e:90091)}
(for Ref.~4, by S. Komlosi) in which it is mentioned explicitly that
Remark 2.3 of Ref.~4 is false by giving a counterexample; \emph{Zbl
1094.26008 Noor, Muhammad Aslam On generalized preinvex functions
and monotonicities. (English) [J] JIPAM, J. Inequal. Pure Appl.
Math. 5, No. 4, Paper No. 110, 9 p., electronic only (2004). ISSN
1443-5756} (by J. E. Mart\'{\i}nez-Legaz) in which it is mentioned that
all the results in the paper follow from a simple observation;
\emph{Zbl 1096.26006 Noor, Muhammad Aslam; Noor, Khalida Inayat On
strongly generalized preinvex functions. (English) [J] JIPAM, J.
Inequal. Pure Appl. Math. 6, No. 4, Paper No. 102, 8 p., electronic
only (2005). ISSN 1443-5756} (by J. E. Mart\'{\i}nez-Legaz) in which,
besides other remarks, it is mentioned a definition which does not
make sense; \emph{Zbl 1093.26006} (for Ref.~7, by N. Hadjisavvas)
where it is mentioned that \textquotedblleft Many other notions and
properties introduced in this paper can be derived in the same way
from the usual generalized invexity notions that can be found in
other papers in the field. When this is not the case, mistakes occur
frequently\textquotedblright. We also pointed out the unprofessional
way some authors answered questions related to their papers. In
conclusion we consider that there are too many papers related to
invexity, much more that the domain deserves. We consider that the
editors of mathematical journals have to pay much more attention
when accepting to publish such papers, taking into account at least
the lack of criticism in the Invexity Community.


\begin{thebibliography}{9}                                                                                                %
\providecommand{\url}[1]{\texttt{#1}} \providecommand{\urlprefix}{URL }

\bibitem {Yang/Yang/Teo:01}\textsc{Yang, X.M.}, \textsc{Yang, X.Q.}, and
\textsc{Teo, K.L.},  \textit{Characterizations and applications of
prequasi-invex functions}. \newblock Journal of Optimization Theory and
Applications, Vol. 110 pp.  645--668, 2001.

\bibitem {Antczak:03}\textsc{Antczak, T.}, \textit{A class of {$B$}-{$(p,r)$%
}-invex functions and mathematical programming}. \newblock Journal of
Mathematical Analysis and Applications, Vol. 286 pp.  187--206, 2003.

\bibitem {Antczak:05}\textsc{Antczak, T.}, \textit{Mean value in invexity
analysis}. \newblock Nonlinear Analysis. Theory, Methods \& Applications,
Vol.~60 pp.  1473--1484, 2005.

\bibitem {Yang/Yang/Teo:03}\textsc{Yang, X.M.}, \textsc{Yang, X.Q.}, and
\textsc{Teo, K.L.},  \textit{Generalized invexity and generalized invariant
monotonicity}. \newblock Journal of Optimization Theory and Applications, Vol.
117 pp.  607--625, 2003.

\bibitem {Yang/Yang/Teo:05}\textsc{Yang, X.M.}, \textsc{Yang, X.Q.}, and
\textsc{Teo, K.L.},  \textit{Criteria for generalized invex monotonicities}.
\newblock European Journal of Operational Research, Vol. 164 pp. 115--119,  2005.

\bibitem {Antczak:07}\textsc{Antczak, T.}, \textit{New optimality conditions
and duality results of {$G$} type in differentiable mathematical programming}.
\newblock Nonlinear Analysis. Theory, Methods \& Applications, Vol.~66 pp.
1617--1632, 2007.

\bibitem {Noor/Noor:06}\textsc{Noor, M.A.} and \textsc{Noor, K.I.},
\textit{Some characterizations of strongly preinvex functions}.
\newblock Journal of Mathematical Analysis and Applications, Vol. 316 pp.
697--706, 2006.

\bibitem {Fan/Guo:07}\textsc{Fan, L.} and \textsc{Guo, Y.}, \textit{On
strongly {$\alpha$}-preinvex functions}. \newblock Journal of Mathematical
Analysis and Applications, Vol. 330 pp.  1412--1425, 2007.
\end{thebibliography}
\end{document}